\documentclass[10pt]{amsart}
\usepackage{amssymb}
\usepackage{amsthm}
\usepackage{amscd}
\usepackage{amsmath}

\newcommand{\qfd}{\hfill $\Box$}

\newtheorem{thm}{Theorem}[section] 

\newtheorem{lemma}[thm]{Lemma}
\newtheorem{remark}[thm]{Remark}

\numberwithin{equation}{section}  

\allowdisplaybreaks

\begin{document}

\title[Subsums]{On The Critical Number of Finite Groups (II)}

\author{Qinghong Wang}
\address{Center for Combinatorics\\
 Nankai University\\
 Tianjin 300071, P.R. China}
 \email{wqh1208@yahoo.com.cn}

 \author{Yongke Qu}
\address{Center for Combinatorics\\
 Nankai University\\
 Tianjin 300071, P.R. China}
 \email{quyongke@sohu.com}

\thanks{AMS 2010: 11B13, 11P50, 20D15. Key
Words: critical number, nilpotent group, additive basis.  }
\thanks{*Corresponding author: Qinghong Wang, Center for Combinatorics,
 Nankai University, Tianjin 300071, P.R. China. \\
E-mail addresses: wqh1208@yahoo.com.cn (Qinghong Wang).\\
\today}

\begin{abstract}
 Let $G$ be a finite group and $S$  a subset
 of $G\setminus\{0\}$. We call $S$ an \emph{additive basis of $G$} if every element of $G$ can
 be expressed as a sum over a nonempty subset in some order. Let $cr(G)$ be the smallest integer $t$ such that every subset of $G\setminus \{0\}$
 of cardinality $t$ is an additive basis of $G$. In this
 paper, we determine $cr(G)$
 for the following cases: $(i)$ $G$ is a finite nilpotent group; $(ii)$ $G$ is a group of even order which possesses a
 subgroup of index $2$.

\end{abstract}

\maketitle

\section{Introduction and Main Results}

Let $G$ be a finite additively written group(not necessarily
commutative). Let $S=\{a_1,\cdots,a_k\}$ be a subset of
$G\setminus\{0\}$. Define
$\sum(S)=\{a_{i_1}+\cdots+a_{i_l}{\mid}i_1,\cdots,i_l\mbox{ are
distinct and }1\leq l\leq k\}$, and for any $1\leq r \leq k$, define
$\sum_r(S)=\{a_{i_1}+\cdots+a_{i_r}{\mid}i_1,\cdots,i_r \mbox{ are
distinct} \}$. We call $S$ an \emph{additive basis of $G$} if $\sum
(S)=G$. The critical number $cr(G)$ of $G$ is the smallest integer
$t$ such that every subset $S$ of $G\setminus\{0\}$ with
${\mid}S{\mid}\geq t$ forms an additive basis of $G$.

Let $Z_n$ be the cyclic group of $n$ elements. $cr(G)$ was first
introduced and studied by Erd\H{o}s and Heilbronn in 1964 for
$G=Z_p$ where $p$ is a prime. With many mathematicians' efforts,
after nearly half a century, $cr(G)$ has been determined for all
finite abelian groups recently (see
\cite{Did-Man}\cite{Did}\cite{Dia-Ham}\cite{Gao-Ham}\cite{Gri}\cite{Fre-Gao-Ger}).

However, the problem to determine $cr(G)$ for $G$ non-abelian is
widely open. So far, we only have the following results in this
direction.

\begin{thm}(\cite{Gao}, \cite{Wang}) Let $G$ be a finite group of order $n$ and $p$ be the
smallest prime divisor of $n$. Then, $cr(G)=n/p+p-2$ providing one
of the following conditions holds,

$(i)$   $G$ is nilpotent, $p\geq 149$ and $n\geq120p^2$;

\smallskip
$(ii)$  There exists a subgroup of $G$ with index $p$ and the other
prime divisor of $n$(if exists) is larger than $6p$,  $p\geq 149$
and $n\geq120p^2$;

\smallskip
$(iii)$ $G$ be a non-abelian group of order $pq\geq 10$ where $q$ is
a prime.
\end{thm}

In this paper we shall determine $cr(G)$ for all groups $G$ as stated
in the abstract by showing the  following two results.

\begin{thm} Let $G$ be a finite nilpotent group of odd order and let $p$ be
 the smallest prime dividing ${\mid}G{\mid}$. If $\frac{|G|}{p}$ is a
composite number then $cr(G)={\mid}G{\mid}/p+p-2$.
\end{thm}

\begin{thm}
Let $G$ be a finite non-abelian group of even order $n$ which
possesses a subgroup of index $2$. Then,

$(i)$  if $n=6$ then $cr(G)=cr(S_3)=4,$ where $S_3$ denotes the symmetric group of six elements;

\smallskip
$(ii)$  $cr(G)=n/2$, otherwise.
\end{thm}

\begin{remark}The proofs of Theorem 1.2 and 1.3 will be heavily based  on
the ideas contained in \cite{Gao-Ham} and \cite{Gri} respectively.
\end{remark}

\begin{remark}
From Theorem 1.1, Theorem 1.2, Theorem 1.3, and the fact that $cr(G)$
has been determined for all finite abelian groups we know  that, the
critical number $cr(G)$ also has been determined for all finite
nilpotent groups and all finite groups of even order which possesses
a subgroup of index two. However,  for  finite groups which contains
no subgroup with index $p$ ($p$ is the smallest prime divisor of
the order of $G$), we even can't guess the exact value of $cr(G)$.
\end{remark}

\section{Notation and Preliminary Lemmas}

\begin{lemma}(\cite{MN})  Let $G$ be a finite group. Let $A$ and $B$ be subsets of
$G$ such that $|A|+|B|>|G|$. Then $A+B=G$, where $A+B = \{a+b|a\in
A, b\in B\}$.

\end{lemma}

M. B. Nathanson(\cite{MN}, Lemma 2.1) stated the conclusion of Lemma
2.1 for abelian groups, but the method used there does work for the
nonabelian groups. For convenience, we repeat the proof here.
\medskip
\ \ \\
{\sl Proof of Lemma 2.1}. For $g\in G$, let $g-B=\{g-b:b\in B\}$. Since
\medskip

\ \ \ \ \ \ \ \ \ \ \ \ \ \ \ \ \ \ \ \ \ \ \ \ \ \ \ \ \ \ \ \ \ \ \
${\mid}G{\mid}\geq{\mid}A\cup(g-B){\mid}$
\medskip

\ \ \ \ \ \ \ \ \ \ \ \ \ \ \ \ \ \ \ \ \ \ \ \ \ \ \ \ \ \ \ \ \ \ \ \ \ \ \ \
$={\mid}A{\mid}+{\mid}g-B{\mid}-{\mid}A\cap(g-B){\mid}$
\medskip

\ \ \ \ \ \ \ \ \ \ \ \ \ \ \ \ \ \ \ \ \ \ \ \ \ \ \ \ \ \ \ \ \ \ \ \ \ \ \ \
$={\mid}A{\mid}+{\mid}B{\mid}-{\mid}A\cap(g-B){\mid},$\\
it follows that $${\mid}A\cap(g-B){\mid}\geq{\mid}A{\mid}+{\mid}B{\mid}-{\mid}G{\mid}\geq1,$$ and so there exist a element $a\in A$ and $b\in B$ such that $g=a+b$. This completes the proof of the lemma.\qfd

\begin{lemma}(\cite{Did})  Let $p,q$ be two primes and $G$ be a finite abelian group
of order $pq$. Let $S$ be a subset of $G$ such that $0\notin S$ and
${\mid}S{\mid}=p+q-1.\ \ Then \sum(S)=G.$

\end{lemma}

Let $G$ be a finite group. Let $B\subset G$ and $x\in G$. As usual,
we write $\lambda_{B}(x)={\mid}(B+x)\setminus B{\mid}$. For any
$B,x,$ Olson proved in \cite{Olson68}
\begin{equation}\label{eq1}
\lambda_{B}(x)=\lambda_{B}(-x)
\end{equation}
and
\begin{equation}\label{eq2}
\lambda_{B}(x)=\lambda_{G\setminus B}(x).
\end{equation}

We use the following property which is implicit in \cite{Olson68}:
Let $G$ be a finite group. Let $S$ be a subset of $G$ such that
$0\notin S$. Put $B=\sum(S)$. For every $y\in S$, we have

\medskip

\begin{center}
$\lambda_{B}(y)={\mid}(\sum(S)+y)\setminus \sum(S){\mid}$
\end{center}

\medskip

\ \ \ \ \ \ \ \ \ \ \ \ \ \ \ \ \ \ \ \ \ \ \ \ \ \ \ \ \ \ \ \ \ \ \ $\leq{\mid}(\sum(S)+y)\setminus (\sum(S\setminus y)+y){\mid}$

\medskip

\ \ \ \ \ \ \ \ \ \ \ \ \ \ \ \ \ \ \ \ \ \ \ \ \ \ \ \ \ \ \ \ \ \ \ $={\mid}\sum(S)\setminus \sum(S\setminus y){\mid}$

\medskip

\ \ \ \ \ \ \ \ \ \ \ \ \ \ \ \ \ \ \ \ \ \ \ \ \ \ \ \ \ \ \ \ \ \ \ $={\mid}\sum(S){\mid}-{\mid} \sum(S\setminus y){\mid}$.

\medskip

By above analysis we get the following inequality
\begin{equation}\label{eq3}
{\mid}\sum(S){\mid}\geq{\mid}\sum(S\setminus
y){\mid}+\lambda_{B}(y).
\end{equation}

We also use the following result of Olson.

\begin{lemma}(\cite{Olson75}) Let $G$ be a finite group and let $S$ be a generating
subset of $G$ such that $0\notin S$. Let $B$ be a subset of $G$ such
that ${\mid}B{\mid}\leq{\mid}G{\mid}/2.$ Then there is a $x\in S$
such that

\begin{center}
  $\lambda_{B}(x)\geq min(({\mid}B{\mid}+1)/2,({\mid}S\cup
-S{\mid}+2)/4).$
\end{center}

\end{lemma}

\begin{lemma}
Let $G$ be a finite group of odd order. Let $S$ be a  subset of $G$
such that $S\cap -S=\emptyset$ and ${\mid}S{\mid}\geq3$. Then
${\mid}\sum(S){\mid}\geq 2{\mid}S{\mid}.$
\end{lemma}

Proof.  We proceed by induction on $|S|$. For ${\mid}S{\mid}=3$, set
$S=\{a,b,c\}$. In order to prove ${\mid}\sum(S){\mid}\geq6$, we
distinguish three cases.

{\bf Case 1.} $a+b=c$. We consider the sequence $(a,b,c,a+c,c+b,a+b+c)$.
If $a+c=b$, then $2a=0$, a contradiction. If $c+b=a$, then $2b=0$, a
contradiction. If $c+b=a+c$, then $b=a$, a contradiction. If
$c+b=a+b+c$, then $b=c$, a contradiction. If $b=a+b+c$, then $a+c=0$, a contradiction. By above analysis we have
that $\{a,b,c,a+c,c+b,a+b+c\}$ is a set, then
${\mid}\sum(S){\mid}\geq6$.

\medskip
{\bf Case 2.} $a+c=b$. The proof is similar to Case 1.

\medskip
{\bf Case 3.} $b+c=a$. The proof is similar to Case 1.

\medskip
{\bf Case 4.} $a+b\neq c$, $a+c\neq b$, $b+c\neq a$. Now we have that either $a,b,c,a+b,
a+c,b+c$ or $a,b,c,a+b,a+c,a+b+c$ are pairwise distinct. This proves the lemma for $|S|=3$.

\medskip

 Now assume that the lemma is true for smaller $|S|$. Set $B=\sum(S)$. Applying Lemma 2.3 to $B$ or
$G\setminus B$ and using (2.2), there exists a $y\in S$ such that
$\lambda_{B}(y)\geq2$. By (2.3),
${\mid}B{\mid}\geq{\mid}\sum(S\setminus
y){\mid}+2\geq2{\mid}S{\mid}.$ This completes the proof. \qfd

 \medskip
 Let $X$ be a subset of $G$ with cardinality $k$. Let
$\{x_{i},1\leq i\leq k\}$ be an ordering of $X$. For $0\leq i\leq
k$, set $X_{i}=\{x_{j}{\mid}1\leq j\leq i\}$ and $B_i=\sum(X_i)$.
The ordering $\{x_{1},\cdots,x_{k}\}$ will be called a
\emph{resolving sequence} of $X$ if for all $i$,
$\lambda_{B_{i}}(x_{i})=max\{\lambda_{B_{i}}(x_{j});1\leq j\leq
i\}$. The \emph{critical index} of the resolving sequence is the
maximal integer $t$ such that $X_{t-1}$ generates a proper subgroup
of $G$.

Clearly, every nonempty subset $S$ not containing $0$ admits a
resolving sequence. Moreover, the critical index is $\geq1.$

We shall write $\lambda_{i}=\lambda_{B_{i}}(x_{i})$. By induction we
have using $(2.3)$ for all $1\leq j\leq k$,

\begin{center}
${\mid}\sum(X){\mid}\geq\lambda_{k}+\cdots+\lambda_{j}+{\mid}B_{j-1}{\mid}.$\end{center}

Put $\delta(m)=0$ if $m$ is odd and $=1$ otherwise. If $\sum(X)<n/2$, by Lemma 2.3,
$\lambda_{i}\geq(i+1+\delta(i))/2$ for all $i\geq t$. In particular
for all $s\geq t$, we have

\begin{equation}\label{eq4}
{\mid}\sum(X){\mid}\geq(k+s+3)(k-s+1)/4-1/2+{\mid}B_{s-1}{\mid}.
\end{equation}

\begin{lemma} Let $G$ be a finite group of order 9, and let $A,B$ be two
subsets of $G$.

$(i)$\ \ If  $|A|=3$ and $A$ is zero-sum free then $|\sum(A)|\geq 6.$

$(ii)$  If $|A|=3$ and $0\not \in A$ then $|\sum(A)|\geq 5.$

$(iii)$If $|A|=4$ and $0\not \in A$ then $|\sum(A)|\geq 7.$

$(iv)$ If $|A|=4$ then $|\sum_2(A)|\geq 5.$

$(v)$   If $|A|=4$ and $|B|\geq 2$ then $|A+B|\geq 5$.
\end{lemma}

\proof By the basic knowledge of group theory we know that $G$ is
abelian.

$(i)$ One can find a proof in \cite{GH}.

$(ii)$ Let $A=\{a_1, a_2, a_3\}$. Assume to the contrary that
$$
|\sum(A)|\leq4.
$$
It follows that $|\{a_1,a_2,a_3\}\cap
\{a_1+a_2,a_1+a_3,a_2+a_3\}|\geq 2.$ Without loss of generality we
may assume that $a_1=a_2+a_3$ and $a_2=a_1+a_3$. Therefore,
$a_1+a_2=a_2+a_3+a_1+a_3$. Hence, $2a_3=0$. Thus, $a_3=0$ for
$|G|=9$, a contradiction with $A\subset G\setminus \{0\}$.

$(iii)$ Let $A=\{a_1, a_2,a_3,a_4\}$. Assume to the contrary that
$$
|\sum(A)|\leq 6.
$$
If $A$ is zero-sum free then $a_1+a_2+a_3+a_4\not \in
\sum(\{a_1,a_2,a_3\})$. By $(i)$ we obtain that $|\sum(A)|\geq
|\{a_1+a_2+a_3+a_4\}\cup
\sum(\{a_1,a_2,a_3\})|=1+|\sum(\{a_1,a_2,a_3\})|\geq 7$, a
contradiction. Hence,
$$
0\in \sum(A).
$$
Therefore,
$$
|\sum(A)\setminus \{0\}|\leq 5.
$$
This together with $(i)$ implies that

(*) $A$ contains no zero-sum free sequence of length 3.

By rearranging if necessary we may assume that $a_1+a_2+a_3\neq 0$.
By (*) we may assume that $a_1+a_2=0$ (by rearranging if necessary).
Since $a_1\neq a_2$, either $a_1+a_3+a_4\neq 0$ or $a_2+a_3+a_4\neq
0$. Without loss of generality we assume that $a_1+a_3+a_4\neq 0$.
It follows from (*) and $a_1+a_2=0$ that $a_3+a_4=0$. Now we have
$$
A=\{a_1,-a_1,a_3,-a_3\}.
$$
Since $\{0\}\cup A=\{0,a_1,-a_1,a_3,-a_3\} \subset \sum(A)$, by the
contrary hypothesis we infer that $\{a_1+a_3,-(a_1+a_3)\}\cap A\neq
\emptyset$. By the symmetry of $A$ we may assume that $a_1+a_3\in
A$. Therefore, $a_1+a_3=-a_1$ or $a_1+a_3=-a_3$. Again by the
symmetry of $A$ we may assume that $a_1+a_3=-a_1$. Thus,
$a_3=-2a_1$. Now we have $A=\{a_1,-a_1,2a_1,-2a_1\}$. Since $|G|=9$
and $a_1\neq -2a_1$, it is easy to see that
$0,a_1,-a_1,2a_1,-2a_1,3a_1,-3a_1$ are $7$ distinct elements from
$\sum(A)$, a contradiction.

$(iv)$ Let $A=\{a_1, a_2,a_3,a_4\}$. Assume to the contrary that $|\sum_{2}(A)|\leq 4$.
It follows that $|\{a_1+a_2,a_1+a_3,a_1+a_4\}\cap
\{a_2+a_3,a_2+a_4,a_3+a_4\}|\geq 2.$ By rearranging if necessary we
assume that $a_1+a_2=a_3+a_4$ and $a_1+a_3=a_2+a_4$. Thus,
$a_1+a_2+a_1+a_3=a_3+a_4+a_2+a_4$. It follows that $a_1=a_4$, a
contradiction.

$(v)$ Let $B=\{b_1,b_2\}$. Assume to the contrary that $|A+B|\leq 4$.
It follows that $|A+B|=4$, and $b_1+A=b_2+A=A+B$. Therefore,
$\sum_{a\in A}(b_1+a)=\sum_{a\in A}(b_2+a)$. Hence,
$|A|b_1+\sum_{a\in A}a=|A|b_2+\sum_{a\in A}a$. Thus, $4b_1=4b_2$ and
$b_1=b_2$, a contradiction. \qfd

\begin{lemma}
If\ \ ${\mid}G{\mid}=27$ then $cr(G)=10$.
\end{lemma}
Proof. We only need to check the case that $G$ is non-abelian. Since
$G$ is a nilpotent  group, $G$ possesses a normal subgroup $K$ of
index $3$. Suppose $G/K=\langle1+K\rangle$. Let $x \in 1+K$ and
$T=(K\setminus \{0\})\cup \{x\}$. It is easy to see that
$-1+K\not\subset\sum(T)$. This shows that $cr(G)\geq10$. So it
suffices to prove that $cr(G)\leq10$. Let $S\subset G\setminus\{0\}$
and ${\mid}S{\mid}=10$. We want to show that $\sum (S)=G$.

From the basic knowledge on $p$-groups (see \cite{TH}) we know that
 there exist exactly four distinct maximal
subgroups of $G$ and each is a normal subgroup of  order 9, and $G$ equals
to the union of these maximal subgroups. Since $|S|=10=2\times 4+2$,
there exists a maximal subgroup $H$ of $G$ such that
$$
|S\cap H|\geq 3.
$$
Now we fix $a\in G\setminus H$. Then, $G=H\cup (a+H)\cup (2a+H)$. It
suffices to prove the following inclusions hold simultaneously:
$$
H\subset \sum(S),\ \ a+H\subset \sum(S),\ \ 2a+H\subset \sum(S).
$$
Let $A=(a+H)\cap S$ and $B=(2a+H)\cap S$. Suppose
$$
A=\{a+a_1,\cdots,a+a_r\},\ \ B=\{2a+b_1, \cdots,2a+b_t\},
$$
where $r\geq t\geq 0$, $r+t=10-|S\cap H|$, and $a_i,b_j\in H$.

Since $H$ is a normal subgroup of $G$, we also have that
$$
a+a_i=a_i'+a,\ \ 2a+b_j=b_j'+2a,
$$
where $a_i',b_j'\in H.$

\medskip
 We
distinguish three cases.

 {\bf Case 1.} $|S\cap H|\geq 5$. By Lemma 2.2 we get
 $$
 H=\sum(S\cap H)\subset \sum(S).
 $$
 Since $|S\cap H|\leq |H\setminus \{0\}|=8$, $|S\cap (G\setminus
 H)|\geq 2$. Therefore, $(\sum(S\cap (G\setminus H)))\cap
 (a+H)\neq \emptyset$ and $(\sum(S\cap (G\setminus H)))\cap
 (2a+H)\neq \emptyset$. It follows from $H=\sum(S\cap H)$ that
 $a+H\subset \sum(S)$ and $2a+H\subset \sum(S)$.

\medskip
 {\bf Case 2.} $|S\cap H|=4$. Now we have
 $$
 r+t=6.
 $$
 By Lemma 2.5$(iii)$ we obtain that
 $$
 |\sum(S\cap H)|\geq 7.
 $$

{\bf Subcase 2.1.} $r=t=3$. Note that $(2a+b_i)+(a+a_i)=b_i'+2a+a+a_i=b_i'+3a+a_i=h+a_i$, where $h=b_i'+3a\in H$. By Lemma 2.1 $\sum(S\cap H)+\{h+a_1,h+a_2,h+a_3\}=H$. Therefore $H\subset\sum(S\cap H)+(2a+b_1)+A\subset\sum(S)$. Again by Lemma 2.1 we have that $a+H\subset A+\sum(S\cap H)\subset\sum(S)$ and $2a+H\subset B+\sum(S\cap H)\subset\sum(S)$.

{\bf Subcase 2.2.} $r\geq4$ and $t\geq1$. Similar to Subcase 2.1 we know that $H\subset\sum(S)$ and $a+H\subset\sum(S)$. Note that $\sum_2(A)\supset\{a+a_1+a+a_2,a+a_1+a+a_3,a+a_1+a+a_4\}$. Therefore, $\sum_2(A)=2a+C$ with $C\subset H$ and ${\mid}C{\mid}\geq3$. By Lemma 2.1, we have $2a+H=2a+C+\sum(S\cap H)=\sum_2(A)+\sum(S\cap H)\subset\sum(S)$.

{\bf Subcase 2.3.} $r=6$. Similarly to Subcase 2.2 one can prove that $a+H\subset\sum(S)$ and $2a+H\subset\sum(S)$. Since $\{a+a_1+a+a_2+a+a_3,a+a_1+a+a_2+a+a_4,a+a_1+a+a_2+a+a_5\}\subset\sum_3(A)$, we infer that ${\mid}\sum_3(A){\mid}\geq3$. Note that $\sum_3(A)\subset H$. By Lemma 2.1, we have $H=\sum_3(A)+\sum(S\cap H)\subset\sum(S)$.

\medskip
{\bf Case 3.} ${\mid}S\cap H{\mid}=3$. By Lemma 2.5 we get
$${\mid}\sum(S\cap H){\mid}\geq5.$$In this case we have $$r+t=7.$$

{\bf Subcase 3.1.} $r=4$ and $t=3$. Note that $A+B=\{a_1'+a,a_2'+a,a_3'+a,a_4'+a\}+\{2a+b_1,2a+b_2,2a+b_3\}=\{a_1'+3a,a_2'+3a,a_3'+3a,a_4'+3a\}
+\{b_1,b_2,b_3\}\subset H$. Since $3a\in H$, by Lemma 2.5$(v)$, ${\mid}A+B{\mid}\geq5$. It follows from Lemma 2.1 that $H=(A+B)+\sum(S\cap H)\subset\sum(S)$. Note that $a+a_i+(2a+b_1)+(a+a_j)=a+a_i+(b_1'+2a)+(a+a_j)=a+a_i+(b_1'+3a+a_j)=a+b_1'+3a+a_i+a_j$. Therefore, $a+b_1'+3a+\sum_2\{a_1,a_2,a_3,a_4\}\subset\sum_3(A\cup B)$. By Lemma 2.5$(iv)$, ${\mid}\sum_2\{a_1,a_2,a_3,a_4\}{\mid}\geq5$. It follows from Lemma 2.1 that $a+H=a+b_1'+3a+\sum_2(A)+\sum(S\cap H)\subset\sum_3(A\cup B)+\sum(S\cap H)\subset\sum(S)$. Note that $2a+b_i+(a+a_k)+(2a+b_j)=2a+b_i+(a_k'+a)+(2a+b_j)=2a+a_k'+3a+b_j+b_i=2a+3a+a_k'+b_j+b_i$. Therefore, $2a+3a+\{a_1',a_2',a_3',a_4'\}+\sum_2\{b_1,b_2,b_3\}\subset\sum_3(A\cup B)$. By Lemma 2.5, ${\mid}\{a_1',a_2',a_3',a_4'\}+\sum_2\{b_1,b_2,b_3\}{\mid}\geq5$. Again by Lemma 2.1 we have $2a+H=2a+3a+\{a_1'+a_2'+a_3'+a_4'\}+\sum_2\{b_1,b_2,b_3\}+\sum(S\cap H)\subset\sum(S)$.

{\bf Subcase 3.2.} $r=5$ and $t=2$. Similarly to above one can prove that $H\subset\sum(S)$ and $a+H\subset\sum(S)$. Since $A+(2a+b_1)+(2a+b_2)\subset2a+H$, by Lemma 2.1 we infer that  $2a+H=A+(2a+b_1)+(2a+b_2)+\sum(S\cap H)\subset\sum(S)$.

{\bf Subcase 3.3.} $r=6$ and $t=1$. Similarly to above one can prove that $H\subset\sum(S)$ and $a+H\subset\sum(S)$. Note that $\{a+a_1\}+\{a+a_2,\cdots,a+a_6\}\subset\sum_2(A)$. Therefore, ${\mid}\sum_2(A){\mid}\geq5$. Now by Lemma 2.1 and $\sum_2(A)\subset2a+H$ we obtain that $2a+H=\sum_2(A)+\sum(S\cap H)\subset\sum(S)$.

{\bf Subcase 3.4.} $r=7$ and $t=0$. Similarly to above one can prove $a+H\subset\sum(S)$ and $2a+H\subset\sum(S)$. Note that $(a+a_1)+(a+a_2)+\{(a+a_3),\cdots,(a+a_7)\}\subset\sum_3(A)\subset H$. Therefore, ${\mid}\sum_3(A){\mid}\geq5$. Again by Lemma 2.1, $H=\sum_3(A)+\sum(S\cap H)\subset\sum(S)$.
\qfd

\section{The Proofs of The Main Results}

\ \ \\
{\sl Proof of Theorem 1.2}.

Set ${\mid}G{\mid}=n$. Since $G$ is a nilpotent  group, $G$
possesses a normal subgroup $K$ of index $p$. Suppose $G/K=\langle1+K\rangle$.
Let $B$ be any subset of $p-2$ elements in $1+K$ and $T=(K\setminus
\{0\})\cup B$. It is easy to see that $-1+K\not\subset\sum(T)$.
This shows that $cr(G)\geq n/p+p-2$. So it suffices to prove that
$cr(G)\leq n/p+p-2$.

\medskip
Let $S$ be any subset of $G\setminus \{0\}$ with cardinality
$|S|=n/p+p-2$. We need to show $\sum(S)=G$. We proceed by induction
on the number of prime divisors of $n$ (counted with multiplicity).
 By the hypothesis we know that $n= 27$ or $n\geq 45$. By Lemma 2.6 we may assume that $n\geq 45$.  Set
$k(n)=(n/p+p-2)/2$. We shall write sometimes $k$ instead of $k(n)$.
Clearly we may partition $S=X\cup Y$ so that
${\mid}X{\mid}={\mid}Y{\mid}=k$, $X\cap -X=Y\cap -Y=\emptyset$
and ${\mid}\sum(X){\mid}\leq{\mid}\sum(Y){\mid}$.

The result holds by Lemma 2.1 if ${\mid}\sum(X){\mid}>n/2$.
Suppose the contrary. Since n is odd, we have

\begin{equation}\label{eq1}
{\mid}\sum(X){\mid}\leq(n-1)/2.
\end{equation}

Let $\{x_{i};1\leq i\leq k\}$ be a resolving sequence for $X$ with
critical index $t$.

By Lemma 2.4 and note that $n\geq45$, in a similar way to the proof
of Theorem $3.1$ in \cite{Gao-Ham} we can prove that

\begin{equation}\label{eq4}
t\geq n/p^2+p.
\end{equation}

Let $H$ be the proper subgroup generate by $X_{t-1}$. Let $p'$ be
the smallest prime divisor of $n/p$. By (3.2), ${\mid}H\cap
S{\mid}\geq n/(pp')+p'-1$. If $n/p$ is the product of more
than two primes, then by the induction hypothesis, $\sum(S\cap
H)=H$. If $n/p$ is the product of two primes, then by Theorem 1.1
and Lemma 2.2, $\sum(S\cap H)=H$.

Since ${\mid}H{\mid}>n/(pp')$, we see easily that
$q={\mid}G/H{\mid}$ is a prime. Since $G$ is nilpotent, $H$ is a
normal subgroup of $G$. Clearly ${\mid}S\setminus H{\mid}\geq q-1$.
Let $a_1,\cdots,a_{q-1}$ be distinct elements from $S\setminus H$.
We denote by $\bar{a_i}$ the image of $a_i$ in $G/H$ under the
canonical morphism.

By the Cauchy-Davenport Theorem(cf.\cite{MN}),
$\{0,\bar{a}_{1}\}+\cdots+\{0,\bar{a}_{q-1}\}=G/H$. It follows that
$\sum(a_1,\cdots,a_{q-1})+H=G$. The theorem now follows since
$\sum(S\cap H)=H.$
\qfd

\ \\\
{\sl Proof of Theorem 1.3}.

Since $G$ possesses a subgroup of index $2$, in a similar way to the
proof of Theorem 1.2 we can show that $cr(G)\geq n/2$. So, it
suffices to prove that $cr(G)\leq n/2$. In a similar way to the proof of
Lemma 2.6 we can checked the theorem for $n\leq 14$(one can find the structures of nonabelian groups for the case in \cite{TH}). Now assume that $n\geq16$. Let $S$ be a
subset of $G\setminus\{0\}$ of size $n/2$. Let $T=S\cup\{0\}$.

Now fix a subgroup $H$ of index $2$. Then, for any $g\in G$,
$H+2g=H$, so that $2g\in H$. Also the sets $T$ and $g-T$ cannot be
disjoint, because of their sizes, so $g$ has a representation as
$t_1+t_2$ with $t_i\in T$. If $g\notin H$, since $2g\in H$, it
means that $t_1\neq t_2$ in its representation $g=t_1+t_2$.
Tossing away $0$, if it is one of $t_i$'s, we have express $g$ as
a subset sum in $S$.

So from now on, we assume $g\in H$, and split the proof into three
cases according to $k=:{\mid}T\cap H{\mid}$.

Case 1. $k\geq(n/2)-1$. obviously.

Case 2. $3\leq k\leq(n/2)-2$. Consider the collection of sums $h+j$
with $h\in T\cap H$ and $j\in T\cap(G\setminus H)$. These
$k({\mid}T{\mid}-k)$ sums belong to $G\setminus H$, so some element
$v$ occurs in this collection with multiplicity at least

\begin{center}
$\lceil\frac{k({\mid}T{\mid}-k)}{{\mid}G\setminus
H{\mid}}\rceil=\lceil\frac{k(n/2+1-k)}{n/2}\rceil\geq\lceil\frac{3(n/2-2)}{n/2}\rceil=3.$

\end{center}

In other words, we can write $v=h_i+j_i$, for $i=1,2,3$, such that
the $h_i(resp., j_i)$ are distinct elements of $T\cap H(resp.,
T\cap(G\setminus H))$. Since $g-v\notin H$, and since as above
$T$ and $(g-v)-T$ are not disjoint, we can write $g-v=h+j$ or
$g-v=j'+h'$ with $h,h'\in H$ and $j,j'\in
T\cap(G\setminus H)$. Pick $i$ so that $h_i\neq h$ and $j_i\neq
j$ or $h_i\neq h'$ and $j_i\neq j'$(which is possible since
there are three choice for $i$). Then we have $g=h+j+h_i+j_i$ or
$g=j'+h'+h_i+j_i$, which is a sum of distinct elements of $T$.
Omitting $0$ as one of the terms, if present, gives a subset from
$S$.

Case 3. $k\leq2$. Now $T$ contains $G\setminus H$, with the possible
exception of a single element $r$. Fix $v\in T\setminus H$. The
$\frac{n}{2}(n/2-1)^2$ sums $x_1+x_2+x_3$ with $x_1,x_2\in
G\setminus (H\cup\{r\})$ and $x_3\in G\setminus H$. In
particular, $g-v$ can be represented $(n/2-1)^2$ ways as such a sum.
Exactly $n/2-1$ of these sums have $x_1=x_2$, $n/2-1$ have
$x_2=x_3$, and $n/2-1$ have $x_1=x_3$. Also, $n/2-1$ of these sums
have $x_1=v$, $n/2-1$ sums have $x_2=v$, and $n/2-1$ sums have
$x_3=v$. Similar $n/2-1$ of these sums have $x_3=r$. There exists a
form $g-v=v+v+x_3$. Thus there remain at least

\begin{center}
$(n/2)^2-7(n/2-1)+1=\frac{(n-2)(n-16)}{4}+1>0.$\end{center}

sums $x_1+x_2+x_3$ equaling $g-v$ with distinct $x_i\in G\setminus
H$ not equal either $v$ or $r$. So there exists a subset sum
representation $g=x_1+x_2+x_3+v.$ This completes the proof of the theorem.
\qfd
\\

\ \ \ \ \ \ \ \ \ \  \ \ \ \ \ \ \ \ \ \ \ \ \ \ \ \ \ \ ACKNOWLEDGEMENTS\\

The authors would like to express sincere thanks to the referees for their many helpful suggestions and also wish to thank their advisor W.Gao for his constructive comments.


\begin{thebibliography}{99}
\bibitem{Erd} P.Erd\H{o}s and H.Heilbronn, {\it On the additon of residue class mod p}, Acta Arith. 9 (1964) 149-159.

\bibitem{Olson68} J.E. Olson, {\it An addition theorem modulo p}, J. Combin. Theory 5 (1968) 45-52.

\bibitem{Did-Man} G.T. Diderrich and H.B. Mann, {\it combinatorial problems in finite
Abelian groups, In ``A survey of Combinatorial Theory''
(J.N.Srivastava et al., Erd\"{o}s.)} pp. 95--100. North-Holland,
Amsterdam, 1973.

\bibitem{Olson75} J.E.Olson, {\it Sums of sets of group elements}, Acta Arith. 28 (1975) 147-156.

\bibitem{Did} G.T. Diderrich, {\it An Addition Theorem for Abelian Groups of Order pq}, Journal of Number Theory, 7 (1975), 33--48.


\bibitem{TH} Thomas W. Hungerford, Algebra, 1980.

\bibitem{Dia-Ham} J.A. Dias da Silva and Y.O. Hamidoune, {\it Cyclic
spaces for Grassmann derivatives and additive theory}, Bull. Lond.
Math. Soc. 26 (1994) 140-146.

\bibitem{Gao} W.Gao, {\it A combinatorial problem on finite groups}, Acta Math SINCA. 38 (1995) 395-399.


\bibitem{GH} A. Geroldinger, F. Halter-Koch, Non-Unique
Factorizations. Algebraic, Combinatorial and Analytic Theory, Pure
and Applied Mathematics, vol. 278, Chapman $\&$ Hall/CRC, 2006.


\bibitem{Gao-Ham}
W. Gao and Y.O. Hamidoune, {\it On additive bases}, Acta Arith. 88 (1999), 233--237.

\bibitem{Gri}
J.R. Griggs, {\it Spanning subset sums for finite abelian groups},
Discrete Math. 229 (2001), 89--99.

\bibitem{MN} M. B. Nathanson, Additive Number Theory, GTM 165, Springer, 1996.
\bibitem{Fre-Gao-Ger} Michael Freeze, Weidong Gao and Alfred Geroldinger, {\it The critical number
of finite abelian groups}, J. Number Theory. 129 (2009), 2766-2777.

\bibitem{Wang} Qinghong Wang and Yongke Qu, {\it On the critical number of finite groups of order $pq$}, Submitted.

\end{thebibliography}
\end{document}